\documentclass[11pt,onecolumn]{IEEEtran}

\usepackage[colorlinks,linkcolor=red,anchorcolor=green,citecolor=blue]{hyperref}
\usepackage[noadjust]{cite}
\usepackage{mathrsfs}               
\usepackage{theorem}				
\usepackage{stfloats}
\usepackage{comment}
\usepackage{graphicx}			
\usepackage[cmex10]{amsmath}
\usepackage{amssymb}
\usepackage{epstopdf}
\makeatother
\newcommand{\R}{{\mathbb R}}
\newcommand{\C}{{\mathbb C}}

\newcommand{\be}{\begin{equation}}
\newcommand{\ee}{\end{equation}}
\newcommand{\pfbox}{\hfill\mbox{$\Box$}}
	\newtheorem{definition}{Definition}
    \newtheorem{theorem}{Theorem}
	
	\newtheorem{lem}{Lemma}
	\newtheorem{remark}{Remark}

\begin{document}
%
\title{State Consensus for Discrete-time
Multi-agent Systems over Time-varying Graphs}
%
%
%
\author{Ji-Lie~Zhang,~Xiang Chen$^*$ and~Guoxiang~Gu,~\IEEEmembership{Fellow,~IEEE}
\thanks{Ji-Lie~Zhang is with the School
of Information Science and Technology,
Southwest Jiaotong University, Chengdu,
Sichuan, 610031, P.R. China.
Xiang Chen is with the Department of Electrical and Computer
Engineering, University of Windsor, Windsor, ON, N9B 3P4, Canada.
Guoxiang Gu is with the School of Electrical Engineering and Computer
Science, Louisiana State University, Baton Rouge, LA 70803-5901, USA.

Corresponding author with email address: xchen@uwindsor.ca
}
}

%
%

\markboth{Journal of XXX,~Vol.~??, No.~?, Month~?? Year ?}%
{Shell \MakeLowercase{\textit{et al.}}: Bare Demo of IEEEtran.cls for Journals}
%



\maketitle

\begin{abstract}
We study the state consensus problem for linear
shift-invariant discrete-time
homogeneous multi-agent systems (MASs) over
time-varying graphs. A novel approach
based on the small gain theorem is
proposed to design the consensus control
protocols for both neutrally stable and
neutrally unstable MASs, assuming the
uniformly connected graphs. It is
shown that the state consensus
can be achieved for neutrally stable MASs
under a weak uniform observability condition;
for neutrally unstable MASs, the state consensus
entails a strong uniform observability condition.
Two numerical examples are worked out
to illustrate our consensus results.
\end{abstract}

\begin{IEEEkeywords}
Discrete-time multi-agent systems,
state consensus, time-varying graph
\end{IEEEkeywords}

\IEEEpeerreviewmaketitle

\section{Introduction}

Consensus control has received
great attention in the control
community since the work of
\cite{osfm,jlm}. It is interesting
to observe that most of the published work
on consensus control has been
focused on continuous-time MASs
\cite{ldch,sh1,wsa}.
For discrete-time MASs, the problem
of state consensus becomes harder.
To be specific, the Laplacian
matrix associated with
the graph topology plays a key role for
the state consensus. Indeed
\cite{yx} shows that for undirected graphs, the
state consensus is hinged on the
maximum eigen-ratio of the nonzero
eigenvalues of the corresponding Laplacian matrix.
This ratio condition can be
further improved by a power of two \cite{gml}.
For the state consensus over time-varying
MASs, initial work has been
focused on single-integrator MASs
over switching topologies \cite{jlm,rb}.
See also \cite{h2017,llylw,noot,sh2012,yh}
for more recent work. It is worth to
mentioning that the latest work \cite{h2017}
provides a very interesting new development
on the state consensus for neutrally stable
MASs over switching
directed graphs: Under a minimum
dwell-time condition,
and assuming that the graph is uniformly connected,
and the associated Laplacian matrix
satisfies certain inequality, the
state consensus can be shown to hold true.
This is the most general result for the state
consensus over switching topologies,
in the best of our knowledge. However the
dwell-time condition and neutrally stable
MASs are the two limiting factors.

We are motivated by the development
for consensus control of discrete-time MASs
over time-varying graphs, and study the
same state consensus problem as in the
literature over more general time-varying
graphs than switching topologies.
In addition we allow
the MAS to have neutrally unstable dynamics.
A new approach based on the small gain
theorem is proposed in this
paper to tackle the state
consensus problem, and develop the
consensus control protocols. Assuming the
same directed graphs as in \cite{h2017}
but without the dwell-time constraint,
we will show that the state consensus can
be achieved for the neutrally stable
MAS under a weak uniform observability
condition. For neutrally unstable MASs,
state consensus entails a strong
uniform observability condition.
The effectiveness of our proposed
consensus control protocols
is illustrated by two simulation examples.

Our paper is organized as follows.
Section II provides the
problem formulation for the state consensus
over time-varying graphs.
Section III presents the consensus results
for neutrally stable MASs, while
Section IV presents the consensus result
for neutrally unstable MASs. Our
state consensus results are illustrated by
two numerical examples.
The paper is concluded in Section V. The notation
is fairly standard with $\R/\C$ standing for
the set of real/complex numbers. For
a matrix $M$ with real or complex entries,
its maximum singular value is denoted by
$\overline{\sigma}(M)$. If $M$ is square,
its $i$th eigenvalue is denoted by $\lambda_i(M)$.
A matrix $A\in\R^{n\times n}$ is said to be
Schur stable, if all its eigenvalues lie
strictly inside the unit circle. For
more general square matrix $A$, its Mahler measure
is defined by
\be\label{mm}
    M(A) := \prod_{i=1}^n \max\{1,|\lambda_i(A)|\}.
\ee
Let $T(z)$ be a stable transfer matrix. Its
${\cal H}_\infty$ norm is bounded, and define by
\be\label{hinf}
    \|T\|_{{\cal H}_\infty}
    := \sup \limits_{|z|> 1} \overline{\sigma}[T(z)].
\ee
Other notation will be made clear as we proceed.

\section{Problem Formulation}

The homogeneous MAS under our
consideration consists of $N$
discrete-time dynamic
systems described by
\begin{align}\label{xik}
    x_i(k + 1) = Ax_i(k) + Bu_i(k),
    \ \ \ x_i(0) = x_{i0},
\end{align}
where $x_i(k) \in \R^n$
is the state vector,
$u_i(k) \in \R^m$
is the control input,
and $x_{i0}$ is the
initial condition of
the $i$th node agent
system for each 
$i \in {\cal N} := \{1,\cdots,N\}$.
The following assumption is made
without loss of generality \cite{yx,h2017}.

\vspace{2mm}

\noindent {\bf Aassumption} (A)\
All eigenvalues of
$A$ lie on the unit circle, and the
pair $(A,B)$ is reachable.

\vspace{2mm}


For leaderless state consensus, we
consider the distributed control
protocol under state feedback:
\begin{align}\label{uik}
    u_i(k) =-\mu F\sum^N_{j=1}
    a_{i,j}(k)[x_i(k)-x_j(k)], \ \ \
    i \in \cal N,
\end{align}
with $\mu>0$ some constant,
and $F\in \R^{m\times n}$ the
state feedback gain
to be specified later,
assuming the accessability
of all the state vectors.
The double array sequence
$\{a_{i,j}(k)\}$ represents
the edge gains of the time-varying
feedback graph, denoted by
${\cal G}(k)$, with
$a_{i,j}(k)\geq 0$ being the
edge gain from the $j$th node agent
system to the $i$th node agent system.
By convention, $a_{i,i}(k)=0$
for all $i \in \cal N$.

Let $x(k) = \mbox{vec}\{x_1(k),
\ldots,x_N(k)\}$ be the global
state vector. The feedback MAS
described by (\ref{xik}) and (\ref{uik})
admits the state space description
\be\label{globalxk}
    x(k + 1) = [I_N \otimes
    A-\mu{\cal L}(k)\otimes(BF)]x(k)
\ee
with ${\cal L}(k)$ being the
Laplacian matrix associated with
the time-varying digraph
${\cal G}(k)$. The state consensus
control requires design of the
state feedback gain $F$ to achieve
\begin{align}\label{consensus}
    \lim_{k\rightarrow\infty}
    [x_i(k) -x_j(k)] = 0\ \forall i,
    j \in \cal N.
\end{align}

By convention, eigenvalues of ${\cal L}(k)$
are arranged in the ascending order according
to their absolute values, i.e.,
\[
    0 = \lambda_1\{{\cal L}(k)\}
    \leq |\lambda_2\{{\cal L}(k)\}| \leq
    \cdots\leq |\lambda_N\{{\cal L}(k)\}|.
\]
The average of the Laplacian matrix
${\cal L}(k)$ over time interval $[k, k + T_{\rm c})$
with integer $T_{\rm c} > 0$ is defined by
\begin{align}\label{Lkb}
    \overline{{\cal L}}_{T_{\rm c}}(k)
    :=\frac{1}{T_{\rm c}}
    \sum^{T_{\rm c}-1}_{i=0}{\cal L}(k+i).
\end{align}
The graph $\overline{{\cal
G}}_{T_{\rm c}}(k)$
corresponding to the average
Laplacian matrix $\overline{{\cal
L}}_{T_{\rm c}}(k)$
can be interpreted as
the union graph over time interval
$[k, k+T_{\rm c})$. The notion
of the uniformly
connected graph
is defined next.  


\begin{definition}\label{defc} A time-varying
digraph is uniformly connected,
if there exists a finite integer $T_{\rm c} > 0$
such that the corresponding
average Laplacian
matrix defined in (\ref{Lkb}) satisfies
$|\lambda_2\{\overline{{\cal L}}_{T_{\rm c}}(k)\}|
> 0 \ \forall k \geq 0$.
\end{definition}


The following assumption is made
on the Laplacian matrix,
which is borrowed from \cite{h2017}.

\vspace{2mm}

\noindent
{\bf Assumption} (L)\
The Laplacian matrix ${\cal L}(k)$
associated with time-varying
digraph ${\cal G}(k)$ satisfies
\be\label{Lineq}
    {\cal L}(k)+{\cal L}(k)'\geq
    \overline{\mu}{\cal L}(k)'{\cal L}(k)
    \ \forall k\geq 0
\ee
for 
some $\overline{\mu}>0$, and 
$\mu$ in (\ref{uik})
satisfies $0<\mu<\overline\mu$.


If the time-varying graph ${\cal G}(k)$
is undirected satisfying
$\lambda_N\{{\cal L}(k)\}\leq 2$
for all $k\geq 0$,
then Assumption (L) holds with
$\overline\mu=1$. Note that
$\lambda_N\{{\cal L}(k)\}\leq 2$ holds,
if all diagonal elements
${\cal L}(k)$ are bounded by 1.


\begin{remark}\label{rem1}
Since each row of ${\cal L}(k)$
sums to zero, 
$\lambda_1\{{\cal L}(k)\}=0$
with eigenvector
$v_1 =\frac{{\rm 1}_N}{\sqrt{N}}$.
Let $\{v_i\}^N_{i=2}$ be chosen
such that $\{v_i\}^N_{i=1}$ form an orthonormal
basis for $\R^N$. Denote
\begin{align*}
    \widehat{V}=
    \left[\begin{array}{ccc}
    v_2 & \cdots & v_N\end{array}\right]
    \in \R^{N\times(N-1)}, \quad
    V = \left[\begin{array}{cc}
    v_1 & \widehat{V}\end{array}\right] \in
    \R^{N\times N}.
\end{align*}
We call $\widehat{{\cal L}}(k)
= \widehat{V}'{\cal L}(k)\widehat{V}
\in \R^{(N-1)\times(N-1)}$ the
{\it dimension reduced Laplacian matrix}.
It follows that
\begin{align}\label{l7}
        V'{\cal L}(k)V = \left[
                  \begin{array}{cc}
                    0 & \ell(k) \\
                    0 & \widehat{\cal L}(k)\\
        \end{array}
        \right],
        \ \ \  \ell(k) = v_1'{\cal L}(k)
        \widehat{V}.
\end{align}
Assumption (L) implies that $\ell(k)=0$
$\forall\ k\geq 0$, although
the corresponding digraph
${\cal G}(k)$ may not be a balanced
time-varying graph, as pointed out
and discussed in \cite{h2017}. There also
holds $\widehat{\cal L}(k)+\widehat{\cal L}(k)'
\geq \overline\mu\widehat{\cal L}(k)'
\widehat{\cal L}(k)$ $\forall k\geq 0$.
\pfbox
\end{remark}


Applying the similarity transformation
$S= V'\otimes I_n$ to the global state
vector $x(k)$ yields
\begin{subequations}\label{xk}
    \begin{align}
    \overline{x}(k + 1)&= A\overline{x}(k)
    - [\mu\ell(k) \otimes (BF)]\widehat{x}(k), \label{xka} \\
    \widehat{x}(k + 1) &= [I_{N-1} \otimes A - \mu\widehat{{\cal L}}(k) \otimes (BF)]\widehat{x}(k), \label{xkb}
\end{align}
\end{subequations}
where $\widehat{x}(k)
= (\widehat{V}'\otimes I_n)x(k)$,
and $\frac{1}{\sqrt{N}}
\overline{x}(k)$ is the average
state vector at time $k$ with
$\overline{x}(k)$ given by
\[
    \overline{x}(k) = (v_1'\otimes I_n)
    x(k) = \frac{1}{\sqrt{N}}
    \sum^N_{i=1}x_i(k)\ \in \R^n.
\]
The next result can be deduced
based on the known results
under Assumption (L) that implies
$\ell(k)=0$ $\forall k\geq 0$.


\begin{lem}\label{lem1}
Under Assumption (L),
feedback MAS (\ref{globalxk})
achieves the state consensus as
defined in (\ref{consensus}), if and only
if the dynamic system described in
(\ref{xkb}) is asymptotically stable.
\end{lem}


\begin{remark}\label{rem2}
The time-varying dynamic system
described in (\ref{xkb})
can be written equivalently as
\begin{align}\label{tvx}
    \widehat{x}(k + 1) =[\widehat{A}-
    \mu\widehat{B}  \widehat{{\cal L}}_m(k)
    \widehat{F}]\widehat{x}(k),
\end{align}
where $\widehat{A} = I_{N-1}\otimes A,
\widehat{B} = I_{N-1}\otimes B,
\widehat{F} = I_{N-1}\otimes F$, and $\widehat{{\cal L}}_m(k) =
\widehat{{\cal L}}(k)\otimes I_m$.
Denote transfer matrix
\[
    \widehat{G}(z) = I_{N-1} \otimes G(z).
    \ \ \ G(z) = F(zI-A)^{-1}B.
\]
Then system (\ref{tvx}) has the
feedback form as illustrated next:

\vspace*{5mm}

\begin{figure}[h]
\begin{center}
\setlength{\unitlength}{3947sp}%
\begingroup\makeatletter\ifx\SetFigFont\undefined%
\gdef\SetFigFont#1#2#3#4#5{%
  \reset@font\fontsize{#1}{#2pt}%
  \fontfamily{#3}\fontseries{#4}\fontshape{#5}%
  \selectfont}%
\fi\endgroup%
\begin{picture}(5424,1599)(2789,-3673)
\thinlines
\put(3751,-2311){\circle{336}}
\put(4501,-2611){\framebox(900,600){\large $\widehat{G}(z)$}}
\put(3901,-2311){\vector( 1, 0){600}}
\put(4501,-3661){\framebox(900,600){\large $-\mu\widehat{\cal L}_m(k)$}}
\put(4501,-3361){\line(-1, 0){750}}
\put(3751,-3361){\vector( 0, 1){900}}
\put(3976,-2211){\makebox(450,375){\large $\widehat{u}(k)$}}
\put(3001,-2311){\vector( 1, 0){600}}
\put(5901,-2211){\makebox(450,375){\large $\widehat{y}(k)$}}
\put(5401,-2311){\vector( 1, 0){1200}}
\put(6151,-2311){\line( 0,-1){1050}}
\put(6151,-3361){\vector(-1, 0){750}}
\end{picture}%

\vspace*{3mm}

{{\bf Fig. 1} \ System (\ref{tvx})
in the feedback form}
\end{center}
\end{figure}

\vspace*{-3mm}

\noindent The state consensus as formulated
in this section is now equivalent to
the asymptotic stability of the
time-varying feedback system in Fig. 1,
in light of Lemma \ref{lem1}. \pfbox
\end{remark}


Clearly, the observability
of $(F,A)$ and reachability of
$(A,B)$ do not ensure the
feedback stabilizability,
due to the time-varying nature
of $\widehat{\cal L}_m(k)$ in the
feedback path. The following weak
notion is crucial.


\begin{definition}\label{defo}
The pair $\{\widehat{\cal
L}_m(k)\widehat{F},\widehat{A}\}$
is (weakly) uniformly observable,
if there exists an integer $T_{\rm o}>0$
such that 
\be\label{obs}
    O(k_0,T_{\rm o}):=\sum^{k_0+T_{\rm o}}_{k=k_0}
    \widehat{A}'^k\widehat{F}'
    \widehat{\cal L}_m(k)'
    \widehat{\cal L}_m(k)
    \widehat{F}\widehat{A}^k>0
    \ \ \forall k_0\geq 0.
\ee
\end{definition}


Under Assumption (A), all eigenvalues
of $A$ lie on the unit circle. If
$A$ has only semi-simple eigenvalues,
then the open-loop system
represented by $\widehat{G}(z)$
is neutrally stable; Otherwise
the open-loop system
represented by $\widehat{G}(z)$
is neutrally unstable. We examine
two different designs of
the state feedback gain $F$
for achieving the state consensus
that is equivalent to stabilization
of the time-varying feedback system
in Fig. 1 in the next two sections,
respectively.

\section{Neutrally Stable MAS}

Under Assumption (A),
the dynamic system represented by
$G(z)=F(zI-A)^{-1}B$ is neutrally stable,
if all eigenvalues of $A$ are
semi-simple. In this case,
the impulse response of $G(z)$
is bounded at all time indices.
In this section, we restrict
our study to the state feedback
gain $F$ specified in the next lemma.


\begin{lem}\label{lem2}
If under Assumption (A), all
eigenvalues of $A$ are semi-simple,
then there exists $X >0$ such that
\be\label{XF}
    {\rm (a)}\ X = A'XA,
    \qquad {\rm (b)}\ I - B'XB \geq 0.
\ee
If $F = B'XA$, then transfer matrix
\[
    G(z) + \frac{1}{2}I =
    F(zI-A)^{-1}B + \frac{1}{2}I
\]
is positive real (PR) \cite{xh}, i.e.,
\begin{align*}
    \left[G(z)+\frac{1}{2}I\right]^* + \left[G(z)+\frac{1}{2}I\right]
    \geq 0\ 
\end{align*}
for almost all $|z|\geq 1$,
and the pair $(F,A)$ is observable.
Moreover $(A-BF)$ is a
Schur stability matrix.
\end{lem}


Proof: Under Assumption (A),
there exists a nonsingular
matrix $S_a$ such that $S_aAS_a^{-1}$
is block diagonal with
each block in the form of
either $\pm 1$ or
\[
    \left[\begin{array}{cr}
    \cos(\theta) & -\sin(\theta)\\
    \sin(\theta) & \cos(\theta)
    \end{array}\right], \ \ \ \theta\in\R.
\]
It follows that $X=\rho S_a'S_a>0$
satisfies (a) of (\ref{XF}) for
any $\rho>0$. Moreover (b) in (\ref{XF})
can be satisfied by taking an appropriate
$\rho>0$. To show that
$G(z) + \frac{1}{2}I$ is PR, we note that for
all $z\in\C$ satisfying $|z|\geq 1$, there hold
\begin{align*}
    \Phi(z) &:=
    G(z)^{\ast} + G(z) + I \\
    &\geq G(z)^{\ast} + G(z) + B'XB \\
    &= B'(\bar{z}I - A')^{-1}F' +
    F(zI - A)^{-1}B + B'XB \\
    &= B'(\bar{z}I - A')^{-1}[A'X(zI - A)
    + (\bar{z}I - A')XA   +
    (\bar{z}I - A')X(zI - A)]
    (zI - A)^{-1}B  \nonumber \\
    &= B'(\bar{z}I - A')^{-1}(|z|^2X -
    A'XA) (zI - A)^{-1}B \\
    &=(|z|^2 - 1)
    [(zI - A)^{-1}B]^{\ast}
    X(zI - A)^{-1}B \geq 0, \ 
\end{align*}
except at $z$ equal to
eigenvalues of $A$.
It is thus concluded that
$G(z) + \frac{1}{2}I$ is indeed PR.
The observability of $(F,A)$ follows
from the reachability of $(A,B)$ by noting
that both $A$ and $X$ are nonsingular.
The Schur stability of $(A-BF)$ follows from
\[
    X-(A-BF)'X(A-BF) = F'(2I-B'XB)F
    \geq F'F£¬
\]
and the observability of $(F,A)$.
Recall that $X>0$. 
\pfbox


\begin{remark}\label{rem3}
The state feedback gain $F$
as specified in Lemma \ref{lem2}
is basically the same as
in \cite{h2017} except
inequality (b) in (\ref{XF}), which
is crucial to proving the PR property of
$G(z) + \frac{1}{2}I$.
While the observability of
$(F,A)$ and the Schur stability of
$(A-BF)$ are similar to those in \cite{h2017},
the PR property is new, and
plays a key role in proving the state
consensus results without the dwell-time
condition. In addition the state
feedback gain $F$ in Lemma \ref{lem2}
can be extended to tackle the
neutrally unstable $G(z)$
for solving the state consensus
problem in the next section. \pfbox
\end{remark}


To make use of the PR property,
define
\be\label{delm}
    \widehat\Delta_m(k):=
    I-\mu\widehat{\cal L}_m(k),
\ee
and use variable substitution
$\widehat{u}(k)= \widehat{v}(k)
-\widehat{F}\widehat{x}(k)$
in (\ref{tvx}). This results in
the state equation
\be\label{ssk}
    \widehat{x}(k+1) =
    (\widehat{A}-\widehat{B}\widehat{F})
    \widehat{x}(k) + \widehat{B}\widehat{v}(k),
\ee
in light of $\widehat{v}(k)=
\widehat{u}(k)+\widehat{y}(k)
=\widehat\Delta_m(k)\widehat{y}(k)$
and $\widehat{y}(k)=\widehat{F}\widehat{x}(k)$.
Hence the feedback system in Fig. 1
can be converted equivalently to the
feedback system in Fig. 2 where
\begin{align}\label{tf}
    \widehat{T}_F (z) = \widehat{F}(zI-
    \widehat{A}+ \widehat{B}\widehat{F})^{-1} \widehat{B}. 
\end{align}

\vspace*{2mm}

\begin{figure}[h]
\begin{center}
\setlength{\unitlength}{3947sp}%
\begingroup\makeatletter\ifx\SetFigFont\undefined%
\gdef\SetFigFont#1#2#3#4#5{%
  \reset@font\fontsize{#1}{#2pt}%
  \fontfamily{#3}\fontseries{#4}\fontshape{#5}%
  \selectfont}%
\fi\endgroup%
\begin{picture}(5424,1599)(2789,-3673)
\thinlines
\put(3751,-2311){\circle{336}}
\put(4501,-2611){\framebox(900,600){\large $\widehat{T}_F(z)$}}
\put(3901,-2311){\vector( 1, 0){600}}
\put(4501,-3661){\framebox(900,600){\large $\widehat\Delta_m(k)$}}
\put(4501,-3361){\line(-1, 0){750}}
\put(3751,-3361){\vector( 0, 1){900}}
\put(3976,-2211){\makebox(450,375){\large $\widehat{v}(k)$}}
\put(3001,-2311){\vector( 1, 0){600}}
\put(5901,-2211){\makebox(450,375){\large $\widehat{y}(k)$}}
\put(5401,-2311){\vector( 1, 0){1200}}
\put(6151,-2311){\line( 0,-1){1050}}
\put(6151,-3361){\vector(-1, 0){750}}
\end{picture}%

\vspace*{5mm}

{{\bf Fig. 2} \ Equivalent feedback
system to Fig. 1}
\end{center}
\end{figure}

\vspace*{-3mm}

The next result shows that transfer matrix
${\widehat{T}}_F(z)$ 
is bounded real (BR) under the hypotheses of
Lemma \ref{lem2}.

\begin{lem}\label{lem3}
Under Assumption (A) and hypotheses
of Lemma \ref{lem2},
$\widehat{T}_F (z)$ in
(\ref{tf}) satisfies
$\|\widehat{T}_F\|_{{\cal H}_\infty} = 1$.
\end{lem}


Proof: The expression of $\widehat{T}_F (z)$ in (\ref{tf}) shows that
\[
    \widehat{T}_F (z) = I_{N-1}\otimes T_F (z),
     \ \ \ T_F (z) = F(zI - A + BF)^{-1}B.
\]
So $\|\widehat{T}_F\|_{{\cal
H}_\infty}=\|T_F\|_{{\cal H}_\infty}$.
It is a known fact \cite{gq} that
\be\label{hinf}
    M(A)^{\frac{1}{m}} \leq
    \inf_F\|T_F\|_{{\cal H}_{\infty}}
    \leq M(A)=1.
\ee
where $M(A)$ is the Mahler measure
of matrix $A$, defined by
\[
    M(A)=\prod_{i=1}^n\max\{
    |\lambda_i(A)|,1\}.
\]
To show that the infimum
in (\ref{hinf}) is achieved by
those $F$ specified in Lemma \ref{lem2},
we note first that
\begin{align*}
    T_F (z) &= F(zI - A)^{-1}B[I
    + F(zI - A)^{-1}B]^{-1} = G(z)[I + G(z)]^{-1}.
\end{align*}
We note next the equivalence of $\|T_F \|_{{\cal H}_\infty} = 1$ to
\begin{align}\label{tfz}
    T_F(z)^{\ast}T_F (z)\leq I \ \forall |z|\geq1.
\end{align}
The proof of the lemma can then be completed by noting the following equivalence relations:
\begin{align*}
    (\ref{tfz}) \ & \Longleftrightarrow \
    G(z)^{\ast}G(z) \leq [I + G(z)]^{\ast}[I + G(z)]\ \forall |z| \geq 1\\
    \ & \Longleftrightarrow \ \left[G(z) + \frac{1}{2}I\right]^{\ast}+\left[G(z) + \frac{1}{2}I\right] \geq 0
    \ \forall |z|\geq 1.
\end{align*}
The last inequality holds
in light of Lemma \ref{lem2}.
\pfbox

\vspace*{2mm}

The main result of this section
is the following theorem on the
consensus condition.


\begin{theorem}\label{th1}
Under Assumptions (A) and (L),
and hypotheses of
Lemma \ref{lem2}, the feedback
MAS described in (\ref{xk})
achieves the state consensus,
if and only if the pair
$\{\widehat{\cal
L}_m(k)\widehat{F},\widehat{A}\}$
is (weakly) uniformly observable.
\end{theorem}


Proof: Under Assumption (A) and
hypotheses of Lemma \ref{lem2},
$G(z)$ and thus $\widehat{G}(z)$
are both PR. It follows from
Lemma \ref{lem3} that
$\|\widehat{T}_F\|_{{\cal H}_\infty}
= 1$. For proving the
sufficiency, suppose that
$\{\widehat{\cal
L}_m(k)\widehat{F},\widehat{A}\}$
is uniformly observable. Then
an integer $T_{\rm o}>0$ exists
such that inequality (\ref{obs})
holds. We claim that,
in reference to Fig. 2,
Assumption (L) implies that
\be\label{linf}
    \|\widehat\Delta_m(k)\|_{{\cal H}_{\infty}}
    := \sup_{\|\widehat{y}(k)\|_{\ell_2}\neq 0}
    \frac{\|\widehat{v}(k)\|_{\ell_2}}
    {\|\widehat{y}(k)\|_{\ell_2}} <1.
\ee
The above inequality holds, if we
can show that $\widehat{\Delta}_m(k)$
is a strict contraction map
from $\ell^m_2 [k_0,\ k_0 + T_{\rm o})$
to $\ell^m_2 [k_0,\ k_0 + T_{\rm o})$
for all $k_0 \geq 0$. This is indeed true
by the following chain of inequalities:
\begin{align*}
    \sum^{k_0+T_{\rm o}}_{k=k_0}
    \|\widehat{v}(k)\|^2
    =& \sum^{k_0+T_{\rm o}}_{k=k_0}
    \widehat{y}(k)'\widehat{\Delta}_m(k)'
    \widehat{\Delta}_m(k)\widehat{y}(k)\\
    =&\sum^{k_0+T_{\rm o}}_{k=k_0}
    \widehat{y}(k)' [I -
    \mu\widehat{\cal L}_m(k)]'[I -
    \mu\widehat{\cal L}_m(k)]\widehat{y}(k)\\
    \leq & \sum^{k_0+T_{\rm o}}_{k=k_0}
    \widehat{y}(k)' [I - \mu(\overline\mu
    -\mu)\widehat{\cal L}_m(k)'
    \widehat{\cal L}_m(k)]\widehat{y}(k)\\
    =&\sum^{k_0+T_{\rm o}}_{k=k_0}
    \left[\|\widehat{y}(k)\|^2
    - \mu(\overline\mu
    -\mu)\|\widehat{\cal L}_m(k)\widehat{y}(k)\|^2\right]\\
    <&\sum^{k_0+T_{\rm o}}_{k=k_0}
    \|\widehat{y}(k)\|^2 \ \ \forall k_0 \geq 0.
\end{align*}
Recall Assumption (L) that implies inequality
$\widehat{\cal L}(k)+\widehat{\cal L}(k)'
\geq \overline\mu\widehat{\cal
L}(k)'\widehat{\cal L}(k)$ $\forall
k\geq 0$ as discussed in Remark \ref{rem1},
and $\overline\mu>\mu>0$.
The last inequality follows from
the weak uniform observability of
$\{\widehat{\cal
L}_m(k)\widehat{F},\widehat{A}\}$.
Indeed if $\widehat{x}(k_0)\neq 0$
but $\|\widehat{\cal L}_m(k)\widehat{y}(k)\|=0$
$\forall$ $k\in[k_0,\ k_0+T_{\rm o})$, then
\[
    \sum^{k_0+T_{\rm o}-1}_{k=k_0}
    \|\widehat{\cal L}_m(k)\widehat{y}(k)\|^2
    = \widehat{x}(k_0)'O(k_0,T_{\rm o})
    \widehat{x}(k_0) = 0,
\]
with $O(k_0,T_{\rm o})$ the
observability gramian defined
in (\ref{obs}), contradicting
the hypothesis on the
uniform observability. Since
$\|\widehat{T}_F\|_{{\cal H}_{\infty}}
= 1$, and $\|\widehat\Delta_m(k)\|_{{\cal H}_{\infty}}<1$, the feedback system
in Fig. 2 is asymptotically stable,
by the well-known small gain theorem,
and thus the state consensus is achieved.

Conversely, if the (weak)
uniform observability of the pair
$\{\widehat{\cal
L}_m(k)\widehat{F},\widehat{A}\}$ does
not hold. Then there exists
$\widehat{x}(k_0)\neq 0$
and $\|\widehat{\cal L}_m(k)\widehat{y}(k)\|=0$
$\forall$ $k\in[k_0,\ k_0+T_{\rm o})$
for each integer $T_{\rm o}>0$.
In this case, $\widehat{u}(k)=0$
for all $k\geq k_0$, in reference to
Fig. 1. As a result, the asymptotic stability
of the corresponding feedback system
does not hold, implying that the
state consensus is not achieved.
This concludes the proof for the
necessity of the weak uniform
observability of $\{\widehat{\cal
L}_m(k)\widehat{F},\widehat{A}\}$
for achieving the
state consensus. \pfbox


\begin{remark}\label{rem4}
The hypotheses of Theorem \ref{th1}
do not include the uniformly
connected graph, as defined in Definition
\ref{defc}. The reason lies in the
fact that the uniform
observability of $\{\widehat{\cal
L}_m(k)\widehat{F},\widehat{A}\}$
implies the uniformly
connected feedback graphs.
The verification
of the uniform observability
seems to be difficult. However
under Assumption (A) and hypothesis
of Lemma \ref{lem2}, unforced
state responses for $\widehat{x}(k)$
are periodic. Consequently
\[
    \|\widehat{\cal L}_m(k)\widehat{y}(k)\|=0
    \ \forall\ k\in[k_0,\ k_0+T_{\rm o})
\]
requires that $\widehat{\cal L}_m(k)$
and thus ${\cal L}(k)$ change
periodically in harmony with those of
$\widehat{x}(k)$. Hence a simple way to
ensure the (weak) uniform observability condition is
to make sure that the time-varying
graph is not only uniformly connected,
but also change aperiodically, or at
least not in harmony with those of
$\widehat{x}(k)$. In \cite{h2017},
this is achieved by considering
switching graphs with dwell-time
$\tau$ greater than or equal to
$\max\{2,d_0\}$ with $d_0>0$
(referred to as the ``controllability index''
in \cite{h2017})
the smallest integer such that
the controllability matrix
$\left[\begin{array}{cccc}
B & AB & \cdots & A^{d_0-1}B
\end{array}\right]$ has the full rank.
The uniform observability condition holds
for such switching graphs,
in light of its necessity by
Theorem \ref{th1}.
\pfbox
\end{remark}

\vspace*{2mm}

\noindent
{\bf Example 1} \ Simulation studies are
carried out for the feedback
MAS with random initial conditions.
Consider the same example as
that in \cite{h2017} with
$N = 4$ agents without the leader
agent where
\[
    A=\left[
    \begin{array}{rr}
    \cos(\pi/4)   & \sin(\pi/4) \\
    -\sin(\pi/4) & \cos(\pi/4) \\
    \end{array}
    \right],\quad
    B=\left[
    \begin{array}{c}
    0 \\
    1\\
    \end{array}
    \right].
\]
It is noted that $X=I$ satisfies (a) and
(b) in (\ref{XF}). Thus
$F=\left[\begin{array}{cc}
-0.7071  &  0.7071 \end{array}\right]$,
in accordance with Lemma \ref{lem2}.
However $A^4=-I$, implying
that each summation term
of the observability
gramian $O(k_0,T_{\rm o})$
in (\ref{obs}) changes
periodically with the
period equal to 4,
if 
the time-varying graph ${\cal G}(k)$
also changes periodically with
the same period of 4, in harmony
with the unforced
state responses of $\widehat{G}(z)$.
Here we consider
\be\label{Lk}
    {\cal L}(k)=\left\{
    \begin{aligned}
    &{\cal L}_0,\quad k=4\kappa,\\
    &{\cal L}_1,\quad k=4\kappa+1,\\
    &{\cal L}_2,\quad k=4\kappa+2,\\
    &{\cal L}_3,\quad k=4\kappa+3,\\
    \end{aligned}\right.
\ee
for $\kappa \in \{0, 1,\ldots\}$
(the dwell-time $\tau=1$), where
\begin{align*}
    {\cal L}_0=\left[
    \begin{array}{rrrr}
     1 & -1 & 0 & 0 \\
    -1 & 1 & 0 & 0 \\
     0 & 0 & 0 & 0 \\
     0 & 0 & 0 & 0 \\
    \end{array}\right], \ \ \
    {\cal L}_1=\left[
    \begin{array}{rrrr}
     0 & 0 & 0 & 0 \\
     0 & 1 & -1 & 0 \\
     0 & -1 & 1 & 0 \\
     0 & 0 & 0 & 0 \\
    \end{array}
    \right],\\
    {\cal L}_2=\left[
    \begin{array}{rrrr}
     0 & 0 & 0 & 0 \\
     0 & 0 & 0 & 0 \\
     0 & 0 & 1 & -1 \\
     0 & 0 & -1 & 1 \\
    \end{array}\right], \ \ \
    {\cal L}_3=\left[
    \begin{array}{rrrr}
     1 & 0 & 0 & -1 \\
     0 & 0 & 0 & 0 \\
     0 & 0 & 0 & 0 \\
     -1 & 0 & 0 & 1 \\
    \end{array}
    \right].
\end{align*}
Basically in each period, node 1 and
node 2 communicate each other first, followed
by node 2 and node 3, then node 3 and node 4,
and finally node 4 and node 1. Therefore
the time-varying graph is uniformly
connected. In addition Assumption (L)
is satisfied with $\overline\mu=1$
and by taking $\mu=0.5$.
However the state consensus is
not achieved, as shown in Fig. 3.


\begin{figure}[!htb]
\centering
    \includegraphics[scale=0.6]{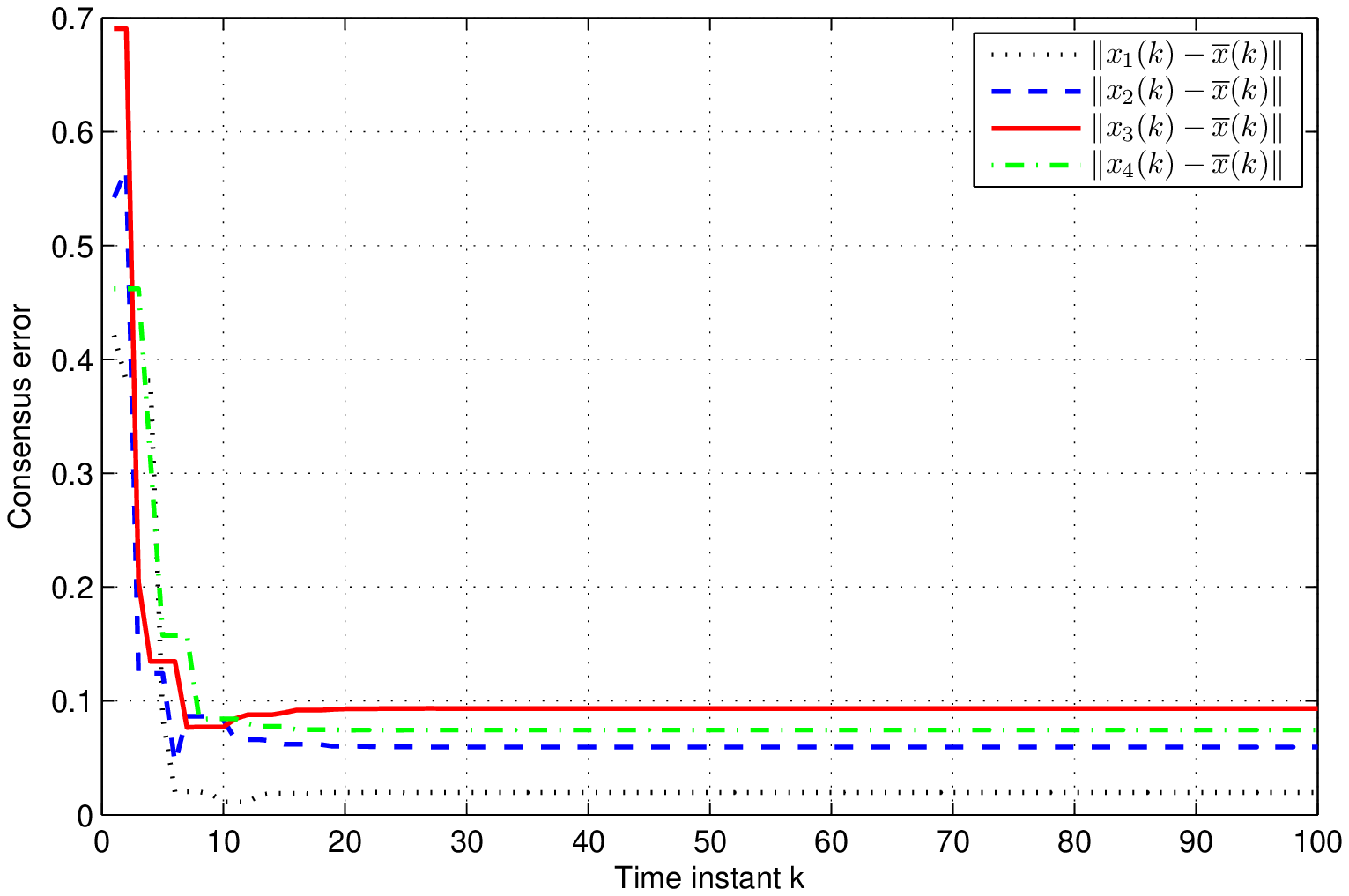}
    \\ {{\bf Fig. 3} \ Consensus error with
    periodic time-varying \\  graph of
    period 4}
\end{figure}


The failure of the state consensus
lies in the fact that the
uniform observability
condition in Theorem \ref{th1}
fails to hold. Indeed, the
observability gramian defined in
(\ref{obs}) over each period is given by
\[
    O(4\kappa,T_{\rm o}=3)=
    0.5\left[
    \begin{array}{rrrrrr}
    1 & 0 & 0 & -1 & 0 & 0 \\
    0 & 1 & -1 & 0 & 0 & 0 \\
    0 & -1 & 2 & 0 & 1 & 0 \\
    -1 & 0 & 0 & 2 & 0 & -1 \\
    0 & 0 & 1 & 0 & 1 & 0 \\
    0 & 0 & 0 & -1 & 0 & 1
    \end{array}\right]
\]
that has rank $4<n(N-1)=6$ for
each integer $\kappa\geq 0$. 
Therefore,
the rank of $O(k_0,T_{\rm o})$ does
not exceed 4 no matter what integer
$T_{\rm o}>0$ is taken.

Suppose that the time-varying
graph ${\cal G}(k)$ involves
periodic switching with period
8 and $\tau=2\geq
\max\{2,d_0=2\}$
for the dwell-time. For instance,
\begin{equation*}
    {\cal L}(k)=\left\{
    \begin{aligned}
    &{\cal L}_0,\quad k=8\kappa,8\kappa+1,\\
    &{\cal L}_1,\quad k=8\kappa+2,8\kappa+3,\\
    &{\cal L}_2,\quad k=8\kappa+4,8\kappa+5,\\
    &{\cal L}_3,\quad k=8\kappa+6,8\kappa+7,\\
    \end{aligned}\right.
\end{equation*}
for $\kappa \in \{0, 1,\ldots\}$. According
to \cite{h2017}, the state consensus
is achieved, validated by the following
simulation plot.


\vspace*{-3mm}

\begin{figure}[!htb]
\centering
    \includegraphics[scale=0.75]{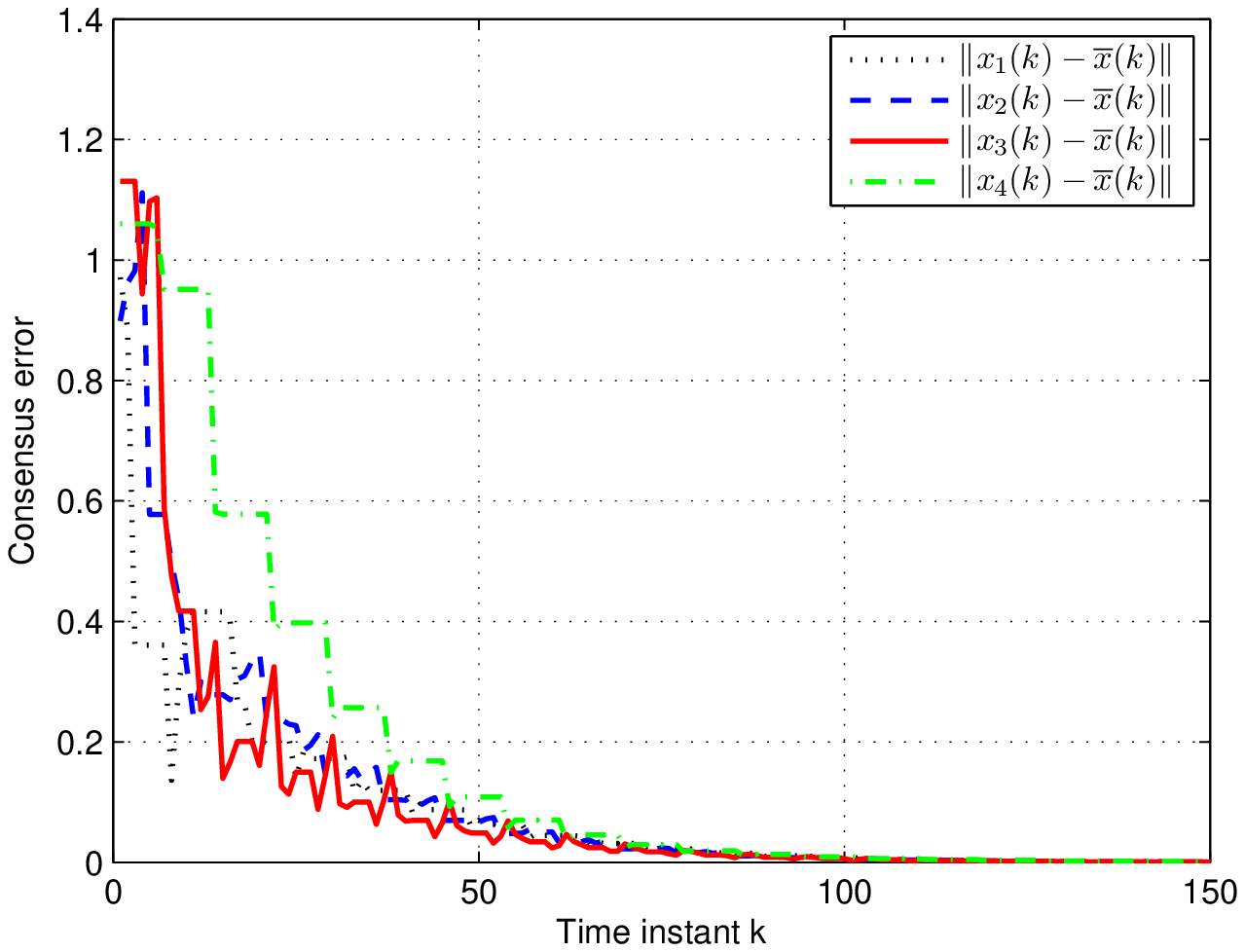}
    \\ {{\bf Fig. 4} \  Consensus error with
    periodic switching graph \\ of
    period 8 and $\tau=2$}
\end{figure}

\noindent In light of Theorem \ref{th1},
the uniform observability condition
holds true, which is verified by
the following observability gramian
computed over the period of 8:
\[
    O(4\kappa,T_{\rm o}=7)=0.25\left[
    \begin{array}{rrrrrr}
    2 & -2 & -2 & 2 & 0 & 0 \\
    -2 & 6 & 2 & -6 & 0 & 0 \\
    -1 & 1 & 8 & 0 & 1 & -1  \\
    1 & -3 & 0 & 8 & -1 & 3 \\
    1 & -1 & 0 & 0 & 7 & 1 \\
    -1 & 3 & 0 & 0 & 1 & 5
  \end{array}
    \right],
\]
for each integer $\kappa\geq 0$.

The uniform observability condition in
Theorem \ref{th1} can also be ensured
by eliminating the periodic
component of the time-varying graph,
in light of Remark \ref{rem4}. This
is indeed true. The next figure
shows the consensus error when the
graph changes randomly among
$\{{\cal L}_i\}_{i=0}^3$ with
equal probability at each
time index $k$.


\vspace*{-3mm}

\begin{figure}[!htb]
\centering
    \includegraphics[scale=0.6]{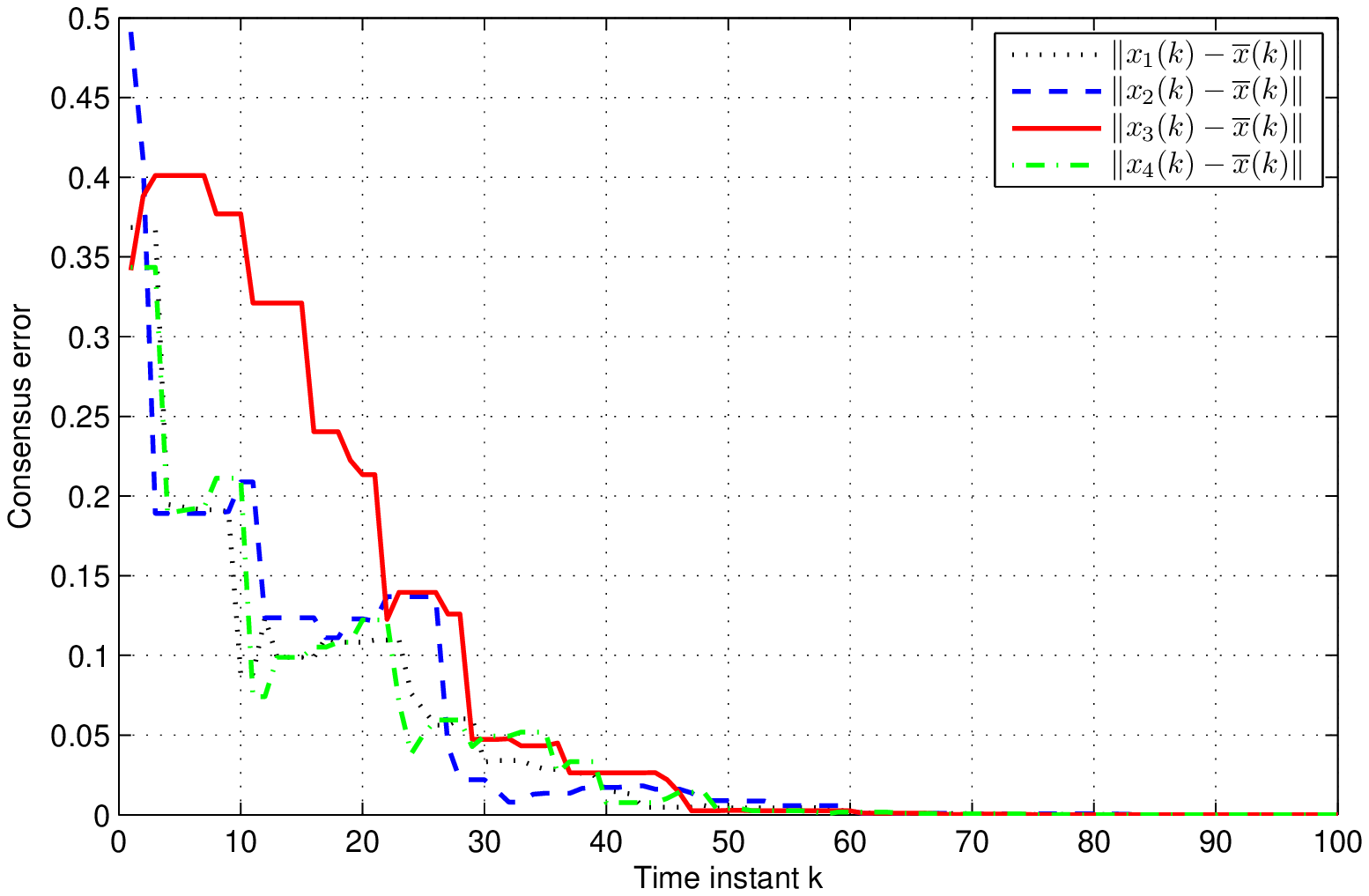}     \\
    {{\bf Fig. 5} \ Consensus error with
     randomly time-varying graph}
\end{figure}

As a concluding remark, we point out
that the uniform observability condition
in Theorem \ref{th1} holds generically.
The exception is when the time-varying
graph changes periodically,
in harmony with the unforced
state responses of $\widehat{G}(z)$,
which needs to be prevented in practice.

\section{Neutrally Unstable MAS}

The uniform observability
condition in Theorem \ref{th1}
requires that the observability
gramian satisfies inequality
(\ref{obs}). However it does not
ensure the exponential convergence
of the consensus error.
For instance, in the case of ${\cal L}(k)
={\cal O}(\frac{1}{k})$ for large
time index $k>0$, then the
exponential convergence fails,
even if the corresponding observability
gramian satisfies inequality (\ref{obs})
for all $k=k_0>0$. 
For this reason we strengthen Definition
\ref{defo} to the following:

\vspace{2mm}

\noindent {\bf Definition 2$'$} \
Recall the observability gramin
in (\ref{obs}).
The pair $\{\widehat{\cal L}_m(k)\widehat
F,\widehat A\}$ is (strongly) uniformly
observable, if there exist an integer
$T_{\rm o}>0$ and a
real number $\varepsilon_{\rm o}>0$
such that
\be\label{obss}
    O(k,T_{\rm o})\geq\varepsilon_{\rm o}I
    \ \ \forall\ k\geq 0.
\ee

\vspace{2mm}

The above strong version of the
uniform observability also
plays a key role
in achieving the state consensus
for neutrally unstable MASs.
Specifically, under Assumption (A)
when not all eigenvalues of
$A$ are semi-simple, the impulse
response of $G(z)=F(zI-A)^{-1}B$
diverges asymptotically. In this
case $X>0$ satisfying (a) and (b)
in (\ref{XF}) of
Lemma \ref{lem2} does not exist.
More importantly, the infimum of
$\|T_F\|_{{\cal H}_{\infty}}$
in (\ref{hinf}) is not achievable
by any stabilizing $F$.
The following lemma provides an extension
to Lemma \ref{lem2} in order to cope with the
neutrally unstable MAS. Its
proof is omitted, since it follows from
the results in \cite{gq}.


\begin{lem}\label{lem4}
Under Assumption (A), for any
given $\gamma>1$, there exists
an $\varepsilon_\gamma>0$
such that algebraic Riccati equation
(ARE)
\be\label{are}
    X = A'X[I + (1-\gamma^{-2})BB'X]^{-1}A
    +\varepsilon_\gamma I,
\ee
admits the stabilizing solution $X>0$,
satisfying $B'XB <\gamma^2I$.
In this case, $\|T_F\|_{{\cal H}_\infty} <
\gamma$ holds with
\be\label{F}
    F = [I + (1 - \gamma^{-2})B'XB]^{-1}B'XA.
\ee
\end{lem}


It is interesting to observe that
with $\gamma=1$ and $\varepsilon_\gamma=0$,
the state feedback gain $F$ in the above
lemma agrees with that in Lemma
\ref{lem2}, which helps to achieve
the state consensus when $A$ has only
semi-simple eigenvalues. The main
result of this section is the following
theorem concerning the state consensus
for the more general case.


\begin{theorem}\label{th2}
Under Assumptions (A) and (L),
let $F$ be designed
as in Lemma \ref{lem4} for
some $\gamma>1$ and
$\varepsilon_{\gamma}>0$,
which is stabilizing and
achieves $\|T_F\|_{{\cal H}_{\infty}}<
\gamma$. If there exists an integer $T_{\rm c}
>0$ and a real number $\epsilon>0$
such that
\be\label{ep}
    \sum^{T_{\rm c}}_{k=0}
    \|\widehat{{\cal L}}_m(k_0+k)\widehat{y}(k_0+k)\|^2
    \geq \epsilon\sum^{T_{\rm c}}_{k=0}
    \|\widehat{y}(k_0+k)\|^2 \ \forall\ k_0\geq 0,
\ee
and $\gamma
\sqrt{1-\epsilon\mu(\overline\mu-\mu)}\leq 1$,
then feedback MAS
(\ref{xk}) achieves the
state consensus.
\end{theorem}


Proof: Assumption (L) being true
implies that the dimension reduced
Laplacian matrix satisfies
inequality (\ref{Lineq}). Note that
$\overline\mu\leq 1$
with equality for undirected time-varying
graphs. In addition $\epsilon>0$ in
Definition $2'$ is necessarily
strictly smaller than 1. Hence
$\sqrt{1-\epsilon\mu(\overline\mu-\mu)}]<1$.
Recall that the state consensus
for the feedback MAS described
in (\ref{globalxk}) is equivalent to
the asymptotic stability of
the feedback system in Fig. 1.
In reference to Fig. 2 that is
equivalent to Fig. 1, there holds
\begin{align*}
    \sum^{k_0+T_{\rm c}}_{k=k_0}
    \|\widehat{v}(k)\|^2
    \leq & \sum^{k_0+T_{\rm c}}_{k=k_0}
    \widehat{y}(k)' [I - \mu(\overline\mu
    -\mu)\widehat{\cal L}_m(k)'
    \widehat{\cal L}_m(k)]\widehat{y}(k)\\
    \leq & \left[1-\epsilon\mu(\overline\mu
    -\mu)\right]\sum^{k_0+T_{\rm c}}_{k=k_0}
    \|\widehat{y}(k)\|^2 \ \ \forall k_0 \geq 0,
\end{align*}
in light of the proof of Theorem \ref{th1}
and inequality (\ref{ep}). 
It follows
that, again in reference to Fig. 2,
\[ 
    \|\widehat\Delta_m\|_{{\cal H}_{\infty}}
    \leq \delta := \sqrt{1-
    \epsilon\mu(\overline\mu-\mu)}.
\] 
Therefore $\|\widehat{T}_F\|_{{\cal
H}_{\infty}}\|\widehat\Delta_m\|_{{\cal H}_{\infty}}<\gamma\delta\leq 1$
by the hypothesis that
$\|T_F\|_{{\cal H}_{\infty}}<\gamma$.
The asymptotic stability of the
feedback system in Fig. 2 follows
from the well-known small gain theorem,
thereby concluding the proof.
\pfbox

\vspace*{2mm}

It is important to observe that
inequality (\ref{ep}) is not
likely to hold under the
notion of the uniform
observability in Definition
\ref{defo}. This is why a stronger
notion in Definition $2'$ is
entailed, which also prevents
${\cal L}(k)$ from decaying to
zero asymptotically, and ensures
the exponential convergence
for the consensus error.

\begin{remark}\label{rem5}
Theoretically, for any $\gamma>1$,
there exists the stabilizing
solution $X>0$ to ARE
(\ref{are}), satisfying $B'XB<\gamma^2I$,
provided that $\varepsilon_{\gamma}>0$
is sufficiently small. However
small $\varepsilon_{\gamma}$ values
lead to low feedback gains of $F$,
and thus long convergence time
for the state consensus errors. In
addition, there is a lower limit for
$\varepsilon_{\gamma}$ in practice,
due to the numerical precision.
Hence $\gamma$ and
$\|\widehat{T}_F\|_{{\cal
H}_{\infty}}$ cannot
be arbitrarily close to one.
Consequently the condition
$\gamma\delta<1$
may not hold in practice, if $\epsilon>0$
is too small, or
if $\varepsilon_{\gamma}>0$ is too small
to exist numerically in Lemma
\ref{lem4} for computing the
stabilizing solution $X>0$
to ARE (\ref{are}) and satisfying
$B'XB<\gamma^2I$.
\pfbox
\end{remark}


\noindent
{\bf Example 2} \ Suppose that
the MAS is the same as in
Example 1, except that
$A$ and $B$ are replaced by
\[
    A=\left[
    \begin{array}{ccc}
    1 & 1 &0 \\
    0 & 1 &1 \\
    0 & 0 &1 \\
    \end{array}
    \right],\quad
    B=\left[
    \begin{array}{r}
    -1 \\
    1\\
    -1\\
    \end{array}
\right].
\]
Clearly the three
eigenvalues of $A$ are all at 1,
and they are not semi-simple.
Hence there does not exist
$F$ such that $\|T_F\|_{{\cal
H}_{\infty}}=1$. By taking $\gamma=1.1$ and $\varepsilon_{\gamma}=10^{-5}$, the
stabilizing solution
$X\geq 0$ to ARE (\ref{are})
satisfying $B'XB<\gamma^2$,
and the corresponding
state feedback gain $F$ in (\ref{F})
can be obtained as follows:
\begin{align*}
    X&=\left[
    \begin{array}{ccc}
     0.0002  &  0.0021  &  0.0103\\
    0.0021  &  0.0304 &  0.1962\\
    0.0103  &  0.1962  &  1.7599\\
    \end{array}
    \right], \\ 
    F&=\left[
    \begin{array}{ccc}
     -0.0068  & -0.1415 &  -1.3985
    \end{array}
    \right].
\end{align*}
We again consider ${\cal L}(k)$ that changes
periodically as in (\ref{Lk}). It can be
verified that the observability gramian
$O(k,T_{\rm o})$ indeed has the full rank
that is 9 with $T_{\rm o}=9$.
In addition $\epsilon>0$
satisfying inequality (\ref{ep})
can be estimated numerically for
which $\gamma\delta<1$
holds true. Recall that $\overline\mu=1$
and $\mu=1/2$ are used in Example 1.
The state consensus errors are shown
 by Fig. 6, and converge
to zero as expected by Theorem \ref{th2}.

\vspace*{-3mm}

\begin{figure}[!htb]
    \centering
    \includegraphics[scale=0.6]{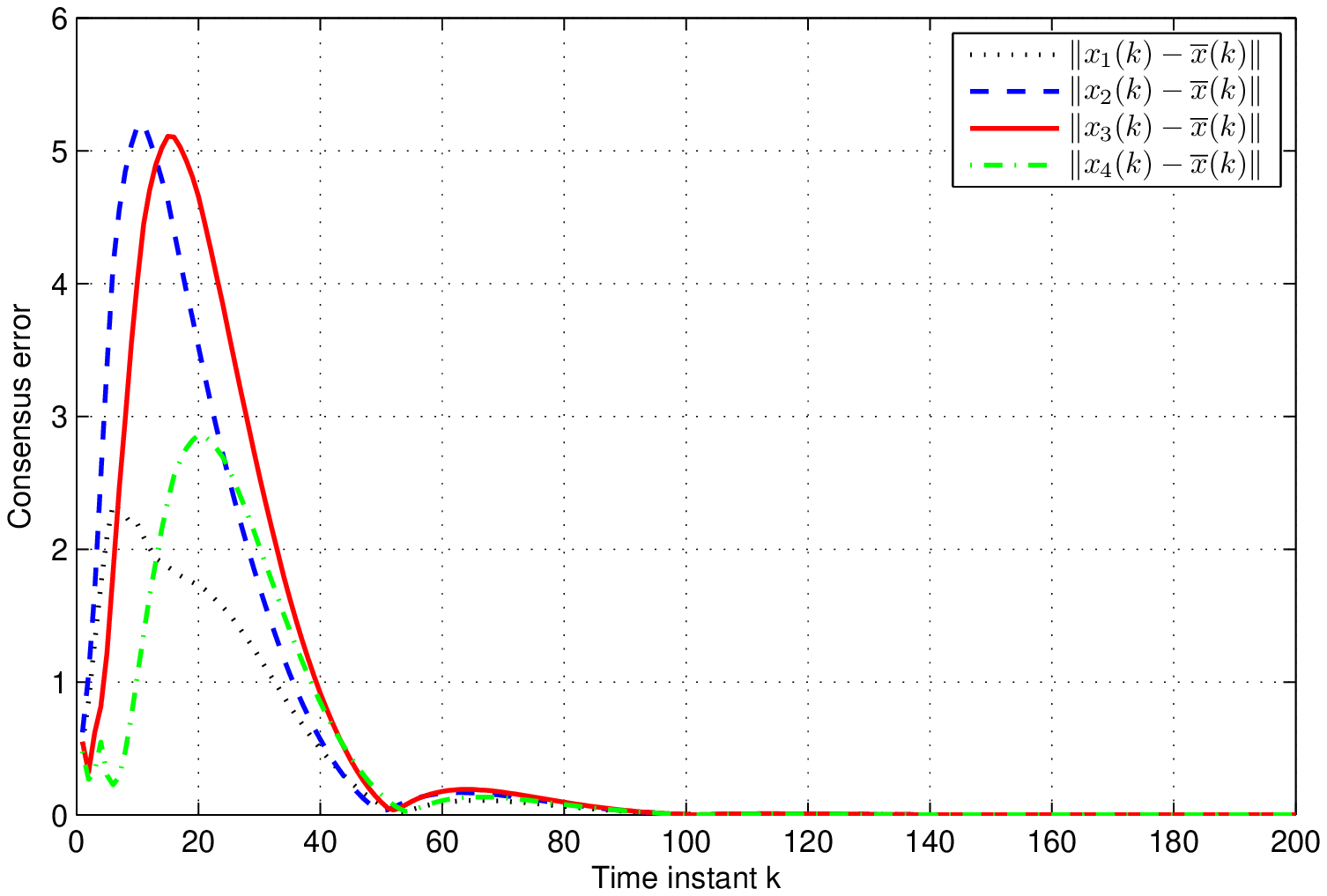}
    \\{{\bf Fig. 6} \ Consensus error for
    neutrally unstable MAS} 
\end{figure}

It is commented that use of Theorem \ref{th2}
to test the small gain condition can be
rather conservative. For instance, the next
simulation study employs the same $F$, but
\be\label{L2k}
    {\cal L}(2\kappa+1)=0, \ \ \
    {\cal L}(2\kappa)={\cal L}_{\kappa\; {\rm mod}\; 4}
\ee
for all positive integers $\{\kappa\}$.
In this case, the
consensus condition in Theorem \ref{th2} fails.
Nevertheless $\|\widehat{T}_F\|_{{\cal H}_{\infty}}\|\widehat\Delta_m(k)\|_{{\cal H}_{\infty}}<1$
holds. The consensus error curves are shown in Fig. 7.


\begin{figure}[!htb]
    \centering
    \includegraphics[scale=0.75]{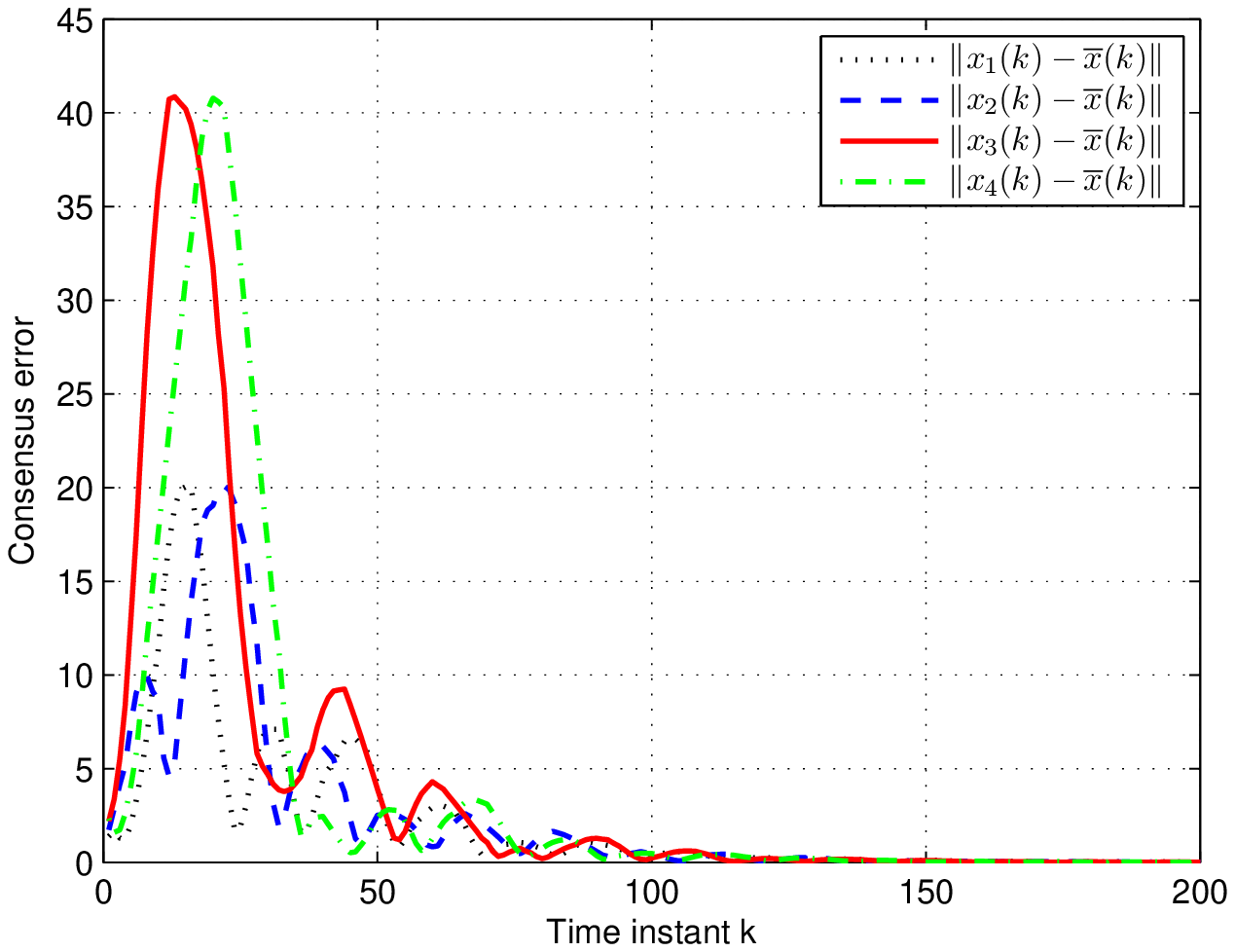}
    \\{{\bf Fig. 7} \ Consensus error for
    neutrally unstable MAS over time-varying
    graph described in (\ref{L2k}) } 
\end{figure}

The above shows that it is
better to estimate $\delta=
\|\widehat\Delta_m(k)\|_{{\cal H}_{\infty}}$
directly. Moreover, the state
consensus takes longer time to achieve than
the previous case, because the time-varying
graph takes twice time intervals to satisfy
the uniform connectivity condition,
and thus the strong uniform observability.

If we use the same $F$ but the time-varying
graph is less frequently connected as follows:
\be\label{L3k}
    {\cal L}(3\kappa+1)={\cal L}(3\kappa+2)=0, \ \ \
    {\cal L}(3\kappa)={\cal L}_{\kappa\; {\rm mod}\; 4}.
\ee
Then even the small gain condition
$\|\widehat{T}_F\|_{{\cal H}_{\infty}}\|\widehat\Delta_m(k)\|_{{\cal H}_{\infty}}<1$
is violated.
The consensus error curves diverge as shown in Fig. 8.


\begin{figure}[!htb]
    \centering
    \includegraphics[scale=0.75]{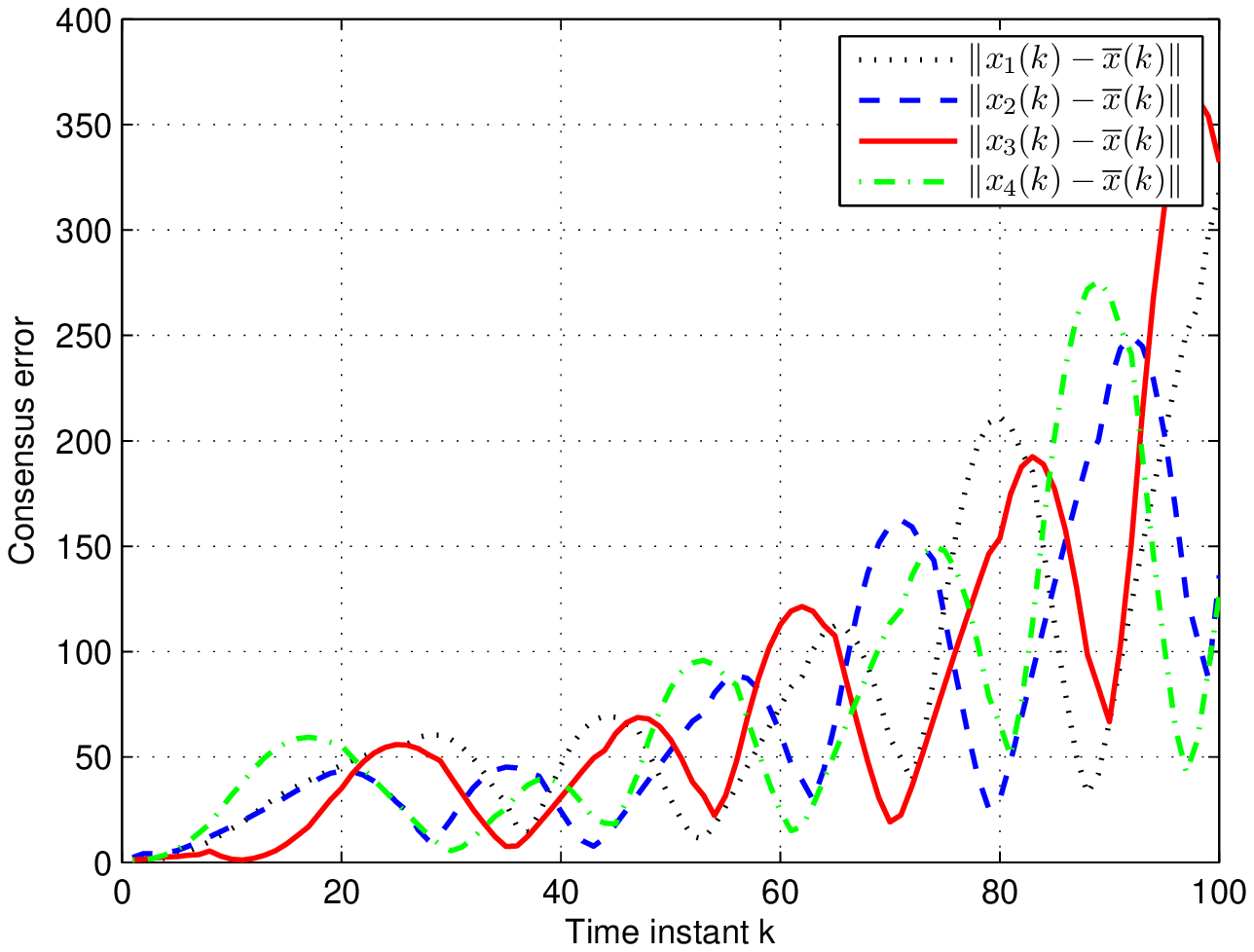}
    \\{{\bf Fig. 8} \ Consensus error for
    neutrally unstable MAS over time-varying
    graph described in (\ref{L3k}) } 
\end{figure}

The above reveals the fact that the neutrally unstable
MAS differs from the neutrally stable MAS. Indeed
the unforced responses diverge for the former, and
remain bounded for the latter. Hence if the control
protocol inputs are frequently absent, then
the consensus errors cannot be pulled back to zero anymore,
no matter what control protocols are used. It suggests
that the time-varying graph needs to be
connected more frequently in order to
achieve the state consensus.

\section{Conclusion}

This paper investigates the problem of
leaderless state consensus for
discrete-time homogeneous MASs
under the state feedback control.
The problem of leader-following
state consensus is not considered
due to the page limit, and essentially
the same mathematical issues. While the
time-varying graphs are the same as
in \cite{h2017} that represents
the latest result on the state consensus
for discrete-time homogeneous MASs
over directed time-varying graphs,
we study both neutrally
stable and neutrally unstable MASs
without the dwell-time constraint
for the time-varying topologies.
A small gain approach is proposed
to develop the control protocols
and to achieve the state consensus.
For the neutrally stable MAS, the PR
property and (weak) uniform observability
condition imply the small gain
condition and ensure the state consensus.
For the neutrally unstable MAS, state consensus
requires a stronger notion
on the uniform observability;
the small gain condition dictates how
often the time-varying graph has to be
connected in order to achieve the
state consensus. Our numerical studies
provide insights to why the
(weak) uniform observability condition
holds generically, and how the
(strong) uniform observability
and small gain condition
can be made true in practice.

\end{document}